\input amstex
\documentstyle{amsppt}
%
\catcode`@=11
\redefine\output@{%
  \def\break{\penalty-\@M}\let\par\endgraf
  \ifodd\pageno\global\hoffset=105pt\else\global\hoffset=8pt\fi  
  \shipout\vbox{%
    \ifplain@
      \let\makeheadline\relax \let\makefootline\relax
    \else
      \iffirstpage@ \global\firstpage@false
        \let\rightheadline\frheadline
        \let\leftheadline\flheadline
      \else
        \ifrunheads@ 
        \else \let\makeheadline\relax
        \fi
      \fi
    \fi
    \makeheadline \pagebody \makefootline}%
  \advancepageno \ifnum\outputpenalty>-\@MM\else\dosupereject\fi
}
\def\Beta{\mathchar"0\hexnumber@\rmfam 42}
\catcode`\@=\active
\nopagenumbers
\chardef\textvolna='176

\chardef\bigalpha='013
\def\negskp{\hskip -2pt}

\chardef\degree="5E
\def\compos{\,\raise 1pt\hbox{$\sssize\circ$} \,}

\def\blue#1{#1}

\catcode`#=11\def\diez{#}\catcode`#=6
\catcode`&=11\catcode`&=4
\catcode`_=11\catcode`_=8
\catcode`\^=11\catcode`\^=7
\catcode`~=11\catcode`~=\active
\def\mycite#1{\cite{\blue{#1}}\immediate\special{ps:
     ShrHPSdict begin /ShrBORDERthickness 0 def}}

\def\mytag#1{%
    \tag#1}
\def\mythetag#1{\thetag{\blue{#1}}\immediate\special{ps:
     ShrHPSdict begin /ShrBORDERthickness 0 def}}
\def\myrefno#1{\no#1}
\def\myhref#1#2{\blue{#2}\immediate\special{ps:
     ShrHPSdict begin /ShrBORDERthickness 0 def}}
\def\myEarXivlink{\myhref{http://arXiv.org}{http:/\negskp/arXiv.org}}

\def\mytheorem#1{\csname proclaim\endcsname{Theorem #1}}
\def\mytheoremwithtitle#1#2{\csname proclaim\endcsname{Theorem #1#2}}
\def\mythetheorem#1{\blue{#1}\immediate\special{ps:
     ShrHPSdict begin /ShrBORDERthickness 0 def}}
\def\mylemma#1{\csname proclaim\endcsname{Lemma #1}}
\def\mylemmawithtitle#1#2{\csname proclaim\endcsname{Lemma #1#2}}

\def\mycorollary#1{\csname proclaim\endcsname{Corollary #1}}

\def\mydefinition#1{\definition{Definition #1}}
\def\mythedefinition#1{\blue{#1}\immediate\special{ps:
     ShrHPSdict begin /ShrBORDERthickness 0 def}}
\def\myconjecture#1{\csname proclaim\endcsname{Conjecture #1}}
\def\myconjecturewithtitle#1#2{\csname proclaim\endcsname{Conjecture #1#2}}

\def\myproblem#1{\csname proclaim\endcsname{Problem #1}}
\def\myproblemwithtitle#1#2{\csname proclaim\endcsname{Problem #1#2}}


\pagewidth{360pt}
\pageheight{606pt}
\vphantom{a}
\vskip -0.8cm
\topmatter
\title
On positive bivariate quartic forms.
\endtitle
\rightheadtext{On positive bivariate quartic forms.}
\author
Ruslan Sharipov
\endauthor
\address Bashkir State University, 32 Zaki Validi street, 450074 Ufa, Russia
\endaddress
\email
\myhref{mailto:r-sharipov\@mail.ru}{r-sharipov\@mail.ru}
\endemail
\abstract
     A bivariate quartic form is a homogeneous bivariate polynomial of 
degree four. A criterion of positivity for such a form is known. In the 
present paper this criterion is reformulated in terms of pseudotensorial 
invariants of the form. 
\endabstract
\subjclassyear{2000}
\subjclass 15A48, 14L24, 11E76\endsubjclass
\endtopmatter
\TagsOnRight
\document

\head
1. Introduction.
\endhead
     Let $a(x^1,x^2)$ be a bivariate quartic form\footnotemark. It is given by 
the following formula:
$$
\hskip -2em
\gathered
a(x^1,x^2)=A_{1111}\,(x^1)^4+4\,A_{1112}\,(x^1)^3\,x^2\,+\\
+\,6\,A_{1122}\,(x^1)^2\,(x^2)^2+4\,A_{1222}\,x^1\,(x^2)^3
+A_{2222}\,(x^2)^4.
\endgathered
\mytag{1.1}
$$
Interpreting $x^1$ and $x^2$ as coordinates of a vector $\bold x$ in
some two-dimensional vector space $V$, one can find that $A_{1111}$, 
$A_{1112}$, $A_{1122}$, $A_{1222}$, $A_{2222}$ are components of a 
symmetric tensor $A$ of the type $(0,4)$. Then the formula \mythetag{1.1} 
is written as 
$$
\hskip -2em
a(x^1,x^2)=\sum^2_{i_1=1}\sum^2_{i_2=1}\sum^2_{i_3=1}\sum^2_{i_4=1}
A_{i_1\,i_2\,i_3\,i_4}\,x^{i_1}\,x^{i_2}\,x^{i_3}\,x^{i_4}.
\mytag{1.2}
$$\par
\footnotetext{\,Upper indices for numerating the variables $x^1$ and $x^2$ 
in \mythetag{1.1} are used according to Einstein's tensorial notation, 
see \mycite{1}.}%
\adjustfootnotemark{-1}%
     The form $a(x^1,x^2)$ in \mythetag{1.1} and in \mythetag{1.2} is 
called positive if it takes positive values for all values of its arguments
$x^1$ and $x^2$ not both zero. The inequalities 
$$
\xalignat 2
&\hskip -2em
A_{1111}>0,
&&A_{2222}>0.
\mytag{1.3}
\endxalignat
$$
are necessary conditions for the positivity of the form \mythetag{1.1}.
However they are not sufficient. Necessary and sufficient conditions
for the positivity of the form \mythetag{1.1} do exist. They are
derived from L.~E.~Dickson's and E.~L.~Rees's results concerning a 
univariate quartic polynomial (see \mycite{2} and \mycite{3}). 
The main goal of the present paper is to express these necessary and 
sufficient conditions in terms of some pseudotensorial invariants 
associated with the tensor $A$ whose components are used in 
\mythetag{1.1}.\par
\head
2. Dickson's and Rees's results.
\endhead
     In \mycite{2} L.~E.~Dickson considers a quartic equation in its 
reduced form with $r\neq 0$:
$$
\hskip -2em
z^4+q\,z^{\kern 0.5pt 2}+r\,z+s=0.
\mytag{2.1}
$$
Among other things, in his book \mycite{2} one can find the following result. 
\mytheorem{2.1} The reduced quartic equation \mythetag{2.1} with real coefficients 
such that $r\neq 0$ has no real roots if and only if its discriminant $D_4>0$ and 
one of the following two conditions is fulfilled:
\roster
\item"1)" \ $4\,s\geqslant q^{\kern 1pt 2}$;
\item"2)" \ $4\,s<q^2$ \ and \ $q\geqslant 0$. 
\endroster
\endproclaim 
\noindent
The discriminant $D_4$ of the reduced quartic equation \mythetag{2.1} is
given by the formula
$$
\hskip -2em
D_4=256\,s^{\kern 1pt 3}-4\,q^{\kern 1pt 3}\,r^{\kern 1pt 2}-27\,r^{\kern 1pt 4}
+16\,q^{\kern 1pt 4}\,s-128\,q^{\kern 1pt 2}\,s^{\kern 1pt 2}
+144\,q\,s\,r^{\kern 1pt 2}.
\mytag{2.2}
$$\par
The case $r=0$ was considered by E.~L.~Rees in \mycite{3}. In this case we have.
\mytheorem{2.2} The reduced quartic equation \mythetag{2.1} with real coefficients 
such that $r=0$ has no real roots if and only if one of the following two conditions
is fulfilled:
\roster
\item"1)" \ $4\,s>q^{\kern 1pt 2}$;
\item"2)" \ $0<4\,s\leqslant q^2$ \ and \ $q>0$. 
\endroster
\endproclaim 
The case $r=0$ is rather simple. In this case, denoting $z^2=y$, one can reduce
the equation \mythetag{2.1} to the following quadratic equation:
$$
\hskip -2em
y^{\kern 1pt 2}+q\,y+s=0.
\mytag{2.3}
$$
The discriminant of the quadratic equation \mythetag{2.3} is given
by the formula 
$$
\hskip -2em
D_2=q^{\kern 1pt 2}-4\,s.
\mytag{2.4}
$$
In the case $r=0$ the discriminant \mythetag{2.2} of the quartic equation 
\mythetag{2.1} reduces to 
$$
\hskip -2em
D_4=16\,s\,(q^{\kern 1pt 2}-4\,s)^2.
\mytag{2.5}
$$
Comparing \mythetag{2.5} with two conditions in Theorem~\mythetheorem{2.2},
we see that the first of them implies $D_4>0$, while the second one implies
$D_4\geqslant 0$.\par
\head
3. Bringing to reduced quartic polynomials.
\endhead 
     Note that the reduced quartic polynomial in the left hand side of
\mythetag{2.1} is positive if and only if the equation \mythetag{2.1}
has no real roots. Therefore bringing the form \mythetag{1.1} to a reduced
quartic polynomial is a method for testing its positivity.\par
     Since $x^1$ and $x^2$ are not both zero in positivity tests, we consider two
cases --- the case $x^1\neq 0$ and the case $x^2\neq 0$. If $x^2\neq 0$, dividing
$a(x^1,x^2)$ by $(x^2)^4>0$ and denoting $t=x^1/x^2$, we derive the following
polynomial from \mythetag{1.1}:
$$
\hskip -2em
P_1(t)=A_{1111}\,t^{\kern 0.5pt 4}+4\,A_{1112}\,t^{\kern  0.5pt 3}
+6\,A_{1122}\,t^{\kern 0.5pt 2}+4\,A_{1222}\,t+A_{2222}.
\mytag{3.1}
$$
If $x^1\neq 0$, dividing $a(x^1,x^2)$ by $(x^1)^4>0$ and denoting $t=x^2/x^1$, 
we get
$$
\hskip -2em
P_2(t)=A_{1111}+4\,A_{1112}\,t+6\,A_{1122}\,t^{\kern  0.5pt 2}
+4\,A_{1222}\,t^{\kern  0.5pt 3}+A_{2222}\,t^{\kern  0.5pt 4}.
\mytag{3.2}
$$
The discriminants $D_{P_{\kern 1.5pt 1}}$  and $D_{P_{\kern 1.5pt 2}}$ 
\pagebreak of the polynomials \mythetag{3.1} and \mythetag{3.2} do coincide. 
They determine the quantity $I_0$ by means of the following formula 
$$
\hskip -2em
I_0=\frac{D_{P_{\kern 1.5pt 1}}}{256}=\frac{D_{P_{\kern 1.5pt 2}}}{256}.
\mytag{3.3}
$$
Writing the formula \mythetag{3.3} explicitly, we get
$$
\hskip -2em
\gathered
I_0=81\,A_{1111}\,(A_{1122})^{\kern 0.5pt 4}\,A_{2222}
-18\,(A_{1111})^{\kern 0.5pt 2}\,(A_{2222})^{\kern 0.5pt 2}
\,(A_{1122})^{\kern 0.5pt 2}\,-\\
-\,27\,(A_{1112})^{\kern 0.5pt 4}\,(A_{2222})^{\kern 0.5pt 2}
-12\,(A_{1111})^{\kern 0.5pt 2}\,A_{1112}\,A_{1222}
\,(A_{2222})^{\kern 0.5pt 2}\,-\\
-\,54\,A_{1111}\,(A_{1122})^{\kern 0.5pt 3}\,(A_{1222})^{\kern 0.5pt 2}
+108\,A_{1111}\,A_{1112}\,(A_{1222})^{\kern 0.5pt 3}\,A_{1122}\,-\\
-\,64\,(A_{1112})^{\kern 0.5pt 3}\,(A_{1222})^{\kern 0.5pt 3}
+54\,A_{1111}\,(A_{2222})^{\kern 0.5pt 2}\,A_{1122}
\,(A_{1112})^{\kern 0.5pt 2}\,+\\
+\,(A_{1111})^{\kern 0.5pt 3}\,(A_{2222})^{\kern 0.5pt 3}
+54\,(A_{1111})^{\kern 0.5pt 2}\,A_{2222}\,A_{1122}
\,(A_{1222})^{\kern 0.5pt 2}\,+\\
+\,36\,A_{1122}^2\,(A_{1112})^{\kern 0.5pt 2}\,(A_{1222})^{\kern 0.5pt 2}
-54\,(A_{1122})^{\kern 0.5pt 3}\,(A_{1112})^{\kern 0.5pt 2}\,A_{2222}\,-\\
-\,27\,(A_{1111})^{\kern 0.5pt 2}\,(A_{1222})^{\kern 0.5pt 4}
-180\,A_{1111}\,A_{1112}\,A_{1222}\,(A_{1122})^{\kern 0.5pt 2}\,A_{2222}\,+\\
+\,108\,(A_{1112})^{\kern 0.5pt 3}\,A_{1222}\,A_{1122}\,A_{2222}
-6\,A_{1111}\,(A_{1112})^{\kern 0.5pt 2}\,(A_{1222})^{\kern 0.5pt 2}\,A_{2222}.
\endgathered
\mytag{3.4}
$$\par
    We can bring the polynomial \mythetag{3.1}
to the reduced form \mythetag{2.1} by substituting 
$$
t=z-\frac{A_{1112}}{A_{1111}}. 
$$
As a result we get the following expressions for the coefficients
in \mythetag{2.1}:
$$
\gather
\hskip -2em
q_1=\frac{6\,A_{1122}}{A_{1111}}-\frac{6\,(A_{1112})^{\kern 0.5pt 2}}
{(A_{1111})^{\kern 0.5pt 2}},
\mytag{3.5}\\
\vspace{1ex}
\hskip -2em
r_1=\frac{4\,A_{1222}}{A_{1111}}
-\frac{12\,A_{1122}\,A_{1112}}{(A_{1111})^{\kern 0.5pt 2}}
+\frac{8\,(A_{1112})^{\kern 0.5pt 3}}{(A_{1111})^{\kern 0.5pt 3}},
\mytag{3.6}\\
\vspace{1ex}
\hskip -2em
s_1=\frac{A_{2222}}{A_{1111}}
-\frac{4\,A_{1222}\,A_{1112}}{(A_{1111})^{\kern 0.5pt 2}}
+\frac{6\,A_{1122}\,(A_{1112})^{\kern 0.5pt 2}}{(A_{1111})^{\kern 0.5pt 3}}
-\frac{3\,(A_{1112})^{\kern 0.5pt 4}}{(A_{1111})^{\kern 0.5pt 4}}.
\mytag{3.7}
\endgather
$$
We can bring the polynomial \mythetag{3.2} to the reduced form 
\mythetag{2.1} by substituting 
$$
t=z-\frac{A_{1222}}{A_{2222}}. 
$$
As a result we get the following expressions for the coefficients
in \mythetag{2.1}:
$$
\gather
\hskip -2em
q_2=\frac{6\,A_{1122}}{A_{2222}}-\frac{6\,(A_{1222})^{\kern 0.5pt 2}}
{(A_{2222})^{\kern 0.5pt 2}},
\mytag{3.8}\\
\vspace{1ex}
\hskip -2em
r_2=\frac{4\,A_{1112}}{A_{2222}}
-\frac{12\,A_{1122}\,A_{1222}}{(A_{2222})^{\kern 0.5pt 2}}
+\frac{8\,(A_{1222})^{\kern 0.5pt 3}}{(A_{2222})^{\kern 0.5pt 3}},
\mytag{3.9}\\
\vspace{1ex}
\hskip -2em
s_2=\frac{A_{1111}}{A_{2222}}
-\frac{4\,A_{1222}\,A_{1112}}{(A_{2222})^{\kern 0.5pt 2}}
+\frac{6\,A_{1122}\,(A_{1222})^{\kern 0.5pt 2}}{(A_{2222})^{\kern 0.5pt 3}}
-\frac{3\,(A_{1222})^{\kern 0.5pt 4}}{(A_{2222})^{\kern 0.5pt 4}}.
\mytag{3.10}
\endgather
$$
\head
4. Positivity criteria for a bivariate quartic form. 
\endhead
     The coefficient $r$ of the reduced polynomial \mythetag{2.1} is used
as a testing parameter in Theorems~\mythetheorem{2.1} and \mythetheorem{2.2}. 
Therefore, relying upon \mythetag{1.3} and using \mythetag{3.6} and 
\mythetag{3.9}, we define the following two testing parameters:
$$
\xalignat 2
&\hskip -2em
I_1=\frac{(A_{1111})^{\kern 0.5pt 3}\,r_1}{4},
&&I_2=\frac{(A_{2222})^{\kern 0.5pt 3}\,r_2}{4}.
\mytag{4.1}
\endxalignat
$$
Here are the explicit expressions for the parameters \mythetag{4.1}:
$$
\align
&\hskip -2em
I_1=(A_{1111})^{\kern 0.5pt 2}\,A_{1222}
-3\,A_{1111}\,A_{1112}\,A_{1122}+2\,(A_{1112})^{\kern 0.5pt 3},
\mytag{4.2}\\
\vspace{1ex}
&\hskip -2em
I_2=(A_{2222})^{\kern 0.5pt 2}\,A_{1112}
-3\,A_{2222}\,A_{1222}\,A_{1122}+2\,(A_{1222})^{\kern 0.5pt 3}.
\mytag{4.3}
\endalign
$$\par
     Another parameter used in Theorems~\mythetheorem{2.1} and 
\mythetheorem{2.2} is $q$. Therefore, relying upon \mythetag{1.3} 
and using \mythetag{3.5} and \mythetag{3.8}, we introduce 
$$
\xalignat 2
&\hskip -2em
I_3=\frac{(A_{1111})^{\kern 0.5pt 2}\,q_1}{6},
&&I_4=\frac{(A_{2222})^{\kern 0.5pt 2}\,q_2}{6}.
\mytag{4.4}
\endxalignat
$$
Here are the explicit expressions for the parameters \mythetag{4.4}:
$$
\align
&\hskip -2em
I_3=A_{1111}\,A_{1122}-(A_{1112})^{\kern 0.5pt 2},
\mytag{4.5}\\
\vspace{1ex}
&\hskip -2em
I_4=A_{2222}\,A_{1122}-(A_{1222})^{\kern 0.5pt 2}.
\mytag{4.6}
\endalign
$$\par
     The third parameter used in Theorems~\mythetheorem{2.1} and 
\mythetheorem{2.2} is $s$. Therefore, relying upon \mythetag{1.3} 
and using \mythetag{3.7} and \mythetag{3.9}, we introduce the following
two parameters: 
$$
\xalignat 2
&\hskip -2em
I_5=(A_{1111})^{\kern 0.5pt 4}\,s_1,
&&I_6=(A_{2222})^{\kern 0.5pt 4}\,s_2.
\mytag{4.7}
\endxalignat
$$
Here are the explicit expressions for the parameters \mythetag{4.7}:
$$
\align
&\hskip -2em
\aligned
I_5=\,\,&6\,A_{1111}\,A_{1122}\,(A_{1112})^{\kern 0.5pt 2}
-3\,(A_{1112})^{\kern 0.5pt 4}-\\
&-\,4\,A_{1111}^{\kern 0.5pt 2}\,A_{1222}\,A_{1112}
+(A_{1111})^{\kern 0.5pt 3}\,A_{2222},
\endaligned
\mytag{4.8}\\
\vspace{1ex}
&\hskip -2em
\aligned
I_6=\,\,&6\,A_{2222}\,A_{1122}\,(A_{1222})^{\kern 0.5pt 2}
-3\,(A_{1222})^{\kern 0.5pt 4}-\\
&-\,4\,A_{2222}^{\kern 0.5pt 2}\,A_{1112}\,A_{1222}
+(A_{2222})^{\kern 0.5pt 3}\,A_{1111}.
\endaligned
\mytag{4.9}
\endalign
$$\par
     Apart from separate entries of $q$ and $s$, Theorems~\mythetheorem{2.1} 
and \mythetheorem{2.2} comprise their combination \mythetag{2.4}. Applying
\mythetag{3.5}, \mythetag{3.7} and \mythetag{3.8}, \mythetag{3.10}
to \mythetag{2.4}, we derive
$$
\allowdisplaybreaks
\gather
\hskip -2em
\aligned
D_{2(1)}&=
\frac{36\,(A_{1122})^{\kern 0.5pt 2}}{(A_{1111})^{\kern 0.5pt 2}}
-\frac{96\,A_{1122}\,(A_{1112})^{\kern 0.5pt 2}}
{(A_{1111})^{\kern 0.5pt 3}}\,+\\
\vspace{1ex}
&+\,\frac{48\,(A_{1112})^{\kern 0.5pt 4}}{(A_{1111})^{\kern 0.5pt 4}}
+\frac{16\,A_{1222}\,A_{1112}}{(A_{1111})^{\kern 0.5pt 2}}
-\frac{4\,A_{2222}}{A_{1111}},
\endaligned
\mytag{4.10}\\
\displaybreak
\hskip -2em
\aligned
D_{2(2)}&=
\frac{36\,(A_{1122})^{\kern 0.5pt 2}}{(A_{2222})^{\kern 0.5pt 2}}
-\frac{96\,A_{1122}\,(A_{1222})^{\kern 0.5pt 2}}
{(A_{2222})^{\kern 0.5pt 3}}\,+\\
\vspace{1ex}
&+\,\frac{48\,(A_{1222})^{\kern 0.5pt 4}}{(A_{2222})^{\kern 0.5pt 4}}
+\frac{16\,A_{1112}\,A_{1222}}{(A_{2222})^{\kern 0.5pt 2}}
-\frac{4\,A_{1111}}{A_{2222}},
\endaligned
\mytag{4.11}
\endgather
$$
Now, relying upon \mythetag{1.3} and using \mythetag{4.10} and 
\mythetag{4.11}, we introduce
$$
\xalignat 2
&\hskip -2em
I_7=\frac{(A_{1111})^{\kern 0.5pt 4}\,D_{2(1)}}{4},
&&I_8=\frac{(A_{2222})^{\kern 0.5pt 4}\,D_{2(2)}}{4}.
\mytag{4.12}
\endxalignat
$$
Here are the explicit expressions for the parameters \mythetag{4.12}:
$$
\align
&\hskip -2em
\aligned
I_7&=9\,(A_{1122})^{\kern 0.5pt 2}\,(A_{1111})^{\kern 0.5pt 2}
-24\,A_{1111}\,A_{1122}\,(A_{1112})^{\kern 0.5pt 2}\,+\\
&+\,12\,(A_{1112})^{\kern 0.5pt 4}
+4\,(A_{1111})^{\kern 0.5pt 2}\,A_{1222}\,A_{1112}
-(A_{1111})^{\kern 0.5pt 3}\,A_{2222},
\endaligned
\mytag{4.13}\\
\vspace{2ex}
&\hskip -2em
\aligned
I_8&=9\,(A_{2222})^{\kern 0.5pt 2}\,(A_{1122})^{\kern 0.5pt 2}
-24\,A_{2222}\,A_{1122}\,(A_{1222})^{\kern 0.5pt 2}\,+\\
&+\,12\,(A_{1222})^{\kern 0.5pt 4}
+4\,(A_{2222})^{\kern 0.5pt 2}\,A_{1112}\,A_{1222}
-(A_{2222})^{\kern 0.5pt 3}\,A_{1111}.
\endaligned
\mytag{4.14}
\endalign
$$
Using the parameters \mythetag{3.4}, \mythetag{4.2}, \mythetag{4.3}, 
\mythetag{4.5}, \mythetag{4.6}, \mythetag{4.6}, \mythetag{4.8}, 
\mythetag{4.9}, \mythetag{4.13}, and \mythetag{4.14}, we can formulate 
the following theorems.
\mytheorem{4.1} A bivariate quartic form \mythetag{1.1} is positive 
if and only if $A_{1111}>0$ and one of the following four conditions 
for its parameters is fulfilled: 
\roster
\item"1)" \ $I_1\neq 0$, \ $I_0>0$, \ $I_7\leqslant 0$;
\item"2)" \ $I_1\neq 0$, \ $I_0>0$, \ $I_7>0$, \ $I_3\geqslant 0$;
\item"3)" \ $I_1=0$, \ $I_7<0$; 
\item"4)" \ $I_1=0$, \ $I_7\geqslant 0$, \ $I_3>0$, \ $I_5>0$. 
\endroster
\endproclaim 
\mytheorem{4.2} A bivariate quartic form \mythetag{1.1} is positive 
if and only if $A_{2222}>0$ and one of the following four conditions 
for its parameters is fulfilled: 
\roster
\item"1)" \ $I_2\neq 0$, \ $I_0>0$, \ $I_8\leqslant 0$;
\item"2)" \ $I_2\neq 0$, \ $I_0>0$, \ $I_8>0$, \ $I_4\geqslant 0$;
\item"3)" \ $I_2=0$, \ $I_8<0$; 
\item"4)" \ $I_2=0$, \ $I_8\geqslant 0$, \ $I_4>0$, \ $I_6>0$. 
\endroster
\endproclaim 
Theorems~\mythetheorem{4.1} and \mythetheorem{4.2} are immediate from
Theorems~\mythetheorem{2.1} and \mythetheorem{2.2}. Each of them is
a criterion of positivity for the bivariate quartic form \mythetag{1.1}.
Hence they are equivalent to each other. However, deriving one of these
theorems from the other seems to be rather difficult.\par
\head
5. Pseudotensors and pseudoscalars. 
\endhead
    Len $V$ be some $n$-dimensional linear vector space. Assume that
$\bold e_1,\,\ldots,\,\bold e_n$ and $\tilde\bold e_1,\,\ldots,
\,\tilde\bold e_n$ are arbitrary two bases in $V$. In this context they 
are usually called the old basis and the new basis respectively (see
\mycite{4}). The bases are related to each other by means of two
mutually inverse matrices square $S$ and $T$: 
$$
\xalignat 2
&\hskip -2em
\tilde\bold e_i=\sum^n_{j=1}S^j_i\,\bold e_j,
&&\bold e_i=\sum^n_{j=1}T^j_i\,\tilde\bold e_j.
\mytag{5.1}
\endxalignat
$$
The matrices $S$ and $T$ are called direct and inverse transition matrices
respectively, while the formulas \mythetag{5.1} are called direct and inverse 
transition formulas. 
\mydefinition{5.1} A pseudotensor of the type $(r,s)$ and 
of the weight $m$ is a geometrical and/or physical object in a linear vector space 
$V$ presented by an array of quantities $F^{i_1\ldots\,i_r}_{j_1\ldots\,j_s}$ 
in each basis $\bold e_1,\,\ldots,\,\bold e_n$ of $V$ and transformed as follows
under any change of basis given by the formulas \mythetag{5.1}: 
$$
F^{i_1\ldots\,i_r}_{j_1\ldots\,j_s}=
(\det T)^m\sum^n\Sb p_1\ldots p_r\\ q_1\ldots q_s\endSb
S^{i_1}_{p_1}\ldots\,S^{i_r}_{p_r}\,\,
T^{q_1}_{j_1}\ldots\,T^{q_s}_{j_s}\,\,
\tilde F^{p_1\ldots\,p_r}_{q_1\ldots\,q_s}.
$$
\enddefinition
     Pseudotensors of the weight $m=0$ are known as tensors (see \mycite{5}). 
\mydefinition{5.2} Pseudotensors of the type $(0,0)$ are called pseudoscalars.
Pseudoscalars of the weight $m=0$ are called scalars.
\enddefinition
     Definition~\mythedefinition{5.1} can be found in \mycite{6}, though 
of cause it was known much prior to \mycite{6} (see \mycite{7} for instance).
This definition is slightly different from that of \mycite{8}.\par
     In our case $\dim V=2$. There is a fundamental pseudotensor $\bold d$ of 
the type $(0,2)$ and of the weight $-1$ in each two-dimensional linear vector
space $V$. Its components are given by the same skew-symmetric matrix 
$$
\hskip -2em
d_{ij}=\Vmatrix\format \r&\quad\l\\ 0 & 1\\
\vspace{1ex}-1 & 0\endVmatrix
$$
in any basis $\bold e_1,\,\bold e_2$ of $V$. The dual object for $\bold d$ is 
given by the same matrix
$$
\hskip -2em
d^{\kern 1pt ij}=\Vmatrix\format \r&\quad\l\\ 0 & 1\\
\vspace{1ex}-1 & 0\endVmatrix
\mytag{5.2}
$$
in any basis $\bold e_1,\,\bold e_2$ of $V$. This dual object is denoted 
by the same symbol $\bold d$ as the initial one. It is a pseudotensor 
of the type $(2,0)$, its weight is equal to $1$. 
\head
6. Pseudotensorial invariants.
\endhead
     The tensor $A$, whose components are used in \mythetag{1.1} and 
\mythetag{1.2}, is a true tensor, i\.\,e\. its weight is zero.
Combining $A$ with the pseudotensor $\bold d$ defined by \mythetag{5.2},
we can compose various pseudotensorial objects. Let's begin with the
following one:
$$
\hskip -2em
B_{i_1\,i_2\,i_3\,i_4}=\sum^2\Sb k_3,\,k_4\\ j_3,\,j_4\endSb
A_{i_1\,i_2\,j_2\,k_1}\,d^{\kern 1pt k_1\,j_1}\,A_{i_3\,i_4\,j_1\,k_2}
\,d^{\kern 1pt k_2\,j_2}
\mytag{6.1}
$$
The formula \mythetag{6.1} defines a pseudotensorial object of the type
$(0,4)$ and of the weight $2$. Its components can be calculated explicitly.
Two of them are associated with the parameters $I_2$ and $I_4$ in
\mythetag{4.5} and \mythetag{4.6}:
$$
\xalignat 2
&\hskip -2em
I_3=-\frac{B_{1111}}{2},
&&I_4=-\frac{B_{2222}}{2}.
\mytag{6.2}
\endxalignat
$$ 
\mydefinition{6.1} Any pseudotensorial object constructed with the use
of $A$ and $\bold d$ is called a pseudotensorial invariant of the quartic
\pagebreak 
form \mythetag{1.1}. 
\enddefinition
So, according to Definition~\mythedefinition{6.1}, the pseudotensor $B$
is a pseudotensorial invariant of the quartic form \mythetag{1.1}. As
for $I_3$ and $I_4$ in \mythetag{6.2}, they are not pseudotensorial 
invariants. They are just certain components of a pseudotensorial 
invariant.\par
     The next pseudotensorial invariant $\hat C$ is constructed with the use 
of $A$, $B$, and $\bold d$ by means of the following formula:
$$
\hskip -2em
\hat C_{i_1\,i_2\,i_3\,i_4\,i_5\,i_6}=\sum^2_{j_4,k_4}
B_{i_1\,i_2\,i_3\,j_4}\,d^{\kern 1pt j_4\,k_4}
\,A_{i_4\,i_5\,i_6\,k_4}.
\mytag{6.3}
$$ 
The formula \mythetag{6.3} defines a pseudotensorial object of the type
$(0,6)$ and of the weight $3$. Its components can be calculated explicitly.
We need only two of them. They are associated with the parameters
$I_1$ and $I_2$ in \mythetag{4.2} and \mythetag{4.3}:
$$
\xalignat 2
&\hskip -2em
I_1=\hat C_{111111},
&&I_2=-\hat C_{222222}.
\mytag{6.4}
\endxalignat
$$\par
     In what follows we need the pseudoscalar object $\beta$ given by
the formula
$$
\hskip -2em
\beta=\sum^2\Sb i_1,\,i_2\\ j_1,\,j_2\endSb
B_{i_1\,i_2\,j_1\,j_2}\,d^{\kern 1pt i_1\,j_1}\,d^{\kern 1pt i_2\,j_2}.
\mytag{6.5}
$$
The pseudoscalar $\beta$ in \mythetag{6.5} can be calculated explicitly:
$$
\hskip -2em
\beta=8\,A_{1112}\,A_{1222}-6\,(A_{1122})^{\kern 0.5pt 2}
-2\,A_{1111}\,A_{2222}.
\mytag{6.6}
$$
Its weight is equal to $4$. Along with $\beta$ in \mythetag{6.6}, we need
the following pseudotensorial object of the type $(0,8)$ and of the weight
$4$:
$$
\hskip -2em
\gathered
D_{i_1\,i_2\,i_3\,i_4\,i_5\,i_6\,i_7\,i_8}\,\,=\!\!\!\sum^2\Sb k_1,k_2,k_3,k_4\\ 
j_1,j_2,j_3,j_4\endSb\!\!\!\!
A_{i_1\,i_2\,j_4\,k_1}\,d^{\kern 1pt k_1\,j_1}
A_{i_3\,i_4\,j_1\,k_2}\,d^{\kern 1pt k_2,j_2}\,\times\\
\vspace{1ex}
\times\,A_{i_5\,i_6\,j_2\,k_3}\,d^{\kern 1pt k_3\,j_3}
\,A_{i_7\,i_8\,j_3\,k_4}\,d^{\kern 1pt k4\,j4}.
\endgathered
\mytag{6.7}
$$
All of the components of the pseudotensorial object \mythetag{6.7} can 
be calculated expli\-citly. We need only two of them:
$$
\aligned
&D_{11111111}=2\,(A_{1112})^{\kern 0.5pt 4}
+2\,(A_{1111})^{\kern 0.5pt 2}\,(A_{1122})^{\kern 0.5pt 2}
-4\,A_{1111}\,A_{1122}\,(A_{1112})^{\kern 0.5pt 2},\\
\vspace{1ex}
&D_{22222222}=2\,(A_{1222})^{\kern 0.5pt 4}
+2\,(A_{2222})^{\kern 0.5pt 2}\,(A_{1122})^{\kern 0.5pt 2}
-4\,A_{2222}\,A_{1122}\,(A_{1222})^{\kern 0.5pt 2}.
\endaligned\quad
\mytag{6.8}
$$
Comparing \mythetag{6.8} with \mythetag{4.8} and \mythetag{4.9} and
using \mythetag{6.6}, we can write
$$
\aligned
I_5=-\frac{3\,D_{11111111}}{2}
-\frac{\beta\,(A_{1111})^{\kern 0.5pt 2}}{2},\\
\vspace{1ex}
I_6=-\frac{3\,D_{22222222}}{2}
-\frac{\beta\,(A_{2222})^{\kern 0.5pt 2}}{2}.
\endaligned
\mytag{6.9}
$$
The formulas \mythetag{6.9} mean that the parameters $I_5$ and $I_6$
are not pseudoscalars. They are just components of the pseudotensorial 
invariant given by the formula 
$$
-\frac{3}{2}\,D_{i_1\,i_2\,i_3\,i_4\,i_5\,i_6\,i_7\,i_8}
-\frac{\beta}{2}\,A_{i_1\,i_2\,i_3\,i_4}\,A_{i_5\,i_6\,i_7\,i_8}.
$$\par
     The parameters $I_7$ and $I_8$ are similar to $I_5$ and $I_6$.
Comparing \mythetag{6.8} with \mythetag{4.13} and \mythetag{4.14} and
taking into account \mythetag{6.6}, we can write
$$
\aligned
I_7=6\,D_{11111111}
+\frac{\beta\,(A_{1111})^{\kern 0.5pt 2}}{2},\\
\vspace{1ex}
I_8=6\,D_{22222222}
+\frac{\beta\,(A_{2222})^{\kern 0.5pt 2}}{2}.
\endaligned
\mytag{6.10}
$$
The formulas \mythetag{6.10} mean that the parameters $I_7$ and $I_8$
are not pseudoscalars. They are just components of the pseudotensorial 
invariant given by the formula 
$$
6\,D_{i_1\,i_2\,i_3\,i_4\,i_5\,i_6\,i_7\,i_8}
+\frac{\beta}{2}\,A_{i_1\,i_2\,i_3\,i_4}\,A_{i_5\,i_6\,i_7\,i_8}.
$$\par
\head
7. Tensorial presentation of the discriminant. 
\endhead
     Let's proceed to the parameter $I_0$ in \mythetag{3.4} which 
was produced from the discriminants $D_{P_{\kern 1.5pt 1}}$ and 
$D_{P_{\kern 1.5pt 2}}$ in \mythetag{3.3}. This parameter is much
more complicated than the previous ones. In dealing with this parameter
we need more sums like \mythetag{6.1} and \mythetag{6.7}. It is convenient
to associate some graphical images with such sums (see Fig\.~7.1) where  
each entry of $A$ corresponds to a node, while each entry of $\bold d$
corresponds to a bond. Each index not used in summation is represented as 
a free bond. Since the tensor $A$ has four indices, each node in
Fig\.~7.1 has exactly four bonds either bound or unbound. For example,
the sum \mythetag{6.1} is presented as dipole with two inner bonds 
and two free bonds at each end.\par
     The pseudoscalar $\beta$ is produced from $B$ according to \mythetag{6.5}.
Graphically the sum \mythetag{6.5} corresponds to binding free bonds of
$B$. Therefore $\beta$ is presented as dipole with four inner bonds in 
Fig\.~7.1. There is an intermediate object 
$$
\hskip -2em
\hat B_{\kern 1pt i1\,i2}=\sum^2_{k_1,j_1} B_{\kern 1pt i1\,j1\,i2\,k1}\,
d^{\kern 1pt k_1\,j_1}.
\mytag{7.1}
$$
It is presented as a dipole with three inner bonds and with one free bond
at each end. One can calculate the components of $\hat B$ in \mythetag{7.1}
explicitly and find that 
$$
\hskip -2em
\hat B_{\kern 1pt i1\,i2}=\frac{\beta}{2}\,d_{\kern 1pt i1\,i2}.
\mytag{7.2}
$$
Due to \mythetag{7.2} the triple dipole $\hat B$ can be replaced with
a bond whenever it enters to a more complicated diagram.\par
    The sum \mythetag{6.7} corresponds to the square $D$ with two free 
bonds at each node. There is one more square shape in Fig\.~7.1. It is
denoted through $\delta$. The shapes with three nodes are presented by 
triangles $C$ and $\gamma$ and by the right angle $\hat C$. The right 
angle $\hat C$ corresponds to the sum \mythetag{6.3}. \par 
     Below we shall derive formulas associated with each shape in Fig\.~7.1. 
Let's begin with the triangular shape $C$. Like the square $D$,
\vadjust{\vskip 560pt\hbox to 0pt{\kern 10pt 
\includegraphics{quartic_form_01.eps}\hss}\vskip 0pt
plus 3pt minus 3pt}the triangle
$C$ has two free bonds at each note. Here is the pseudotensorial 
object associated with $C$:
$$
C_{i_1\,i_2\,i_3\,i_4\,i_5\,i_6}\,\,=\!\!\!\sum^2\Sb k_1,k_2,k_3\\ 
j_1,j_2,j_3\endSb\!\!\!\!
A_{i_1\,i_2\,j_3\,k_1}\,d^{\kern 1pt k_1\,j_1}
A_{i_3\,i_4\,j_1\,k_2}\,d^{\kern 1pt k_2,j_2}
\,A_{i_5\,i_6\,j_2\,k_3}\,d^{\kern 1pt k_3\,j_3}.
\quad
\mytag{7.3}
$$
Its type is $(0,6)$, its weight is $3$. The components of the pseudotensor
\mythetag{7.3} can be calculated explicitly. Using them, we calculate the
pseudoscalar $\gamma$: 
$$
\hskip -2em
\gamma\,\,=\!\!\sum^2\Sb k_1\!,k_2,k_3\\j_1\!,j_2,j3\endSb
C_{j_3\,k_1\,j_1\,k_2\,j_2\,k_3}\,d^{\kern 1pt k_1\,j_1}
\,d^{\kern 1pt k_2,j_2}\,d^{\kern 1pt k_3\,j_3}. 
\mytag{7.4}
$$
The weight of the pseudoscalar \mythetag{7.4} is $6$. Here is the 
explicit formula for $\gamma$:
$$
\aligned
\gamma&=12\,A_{1112}\,A_{1122}\,A_{1222}
+6\,A_{1111}\,A_{1122}\,A_{2222}\,-\\
\vspace{1ex}
&-\,6\,(A_{1122})^{\kern 1pt 3}
-6\,(A_{1112})^{\kern 1pt 2}\,A_{2222}
-6\,A_{1111}\,(A_{1222})^{\kern 1pt 2}. 
\endaligned
$$
\vadjust{\vskip 120pt
\hbox to 0pt{\kern 10pt \includegraphics{quartic_form_02.eps}\hss}
\vskip 5pt}\par
     There is one more triangular shape in Fig\.~7.1. It is denoted $\check C$. 
Here is the formula for the associated pseudotensorial object:
$$
\hskip -2em
\check C_{i_1\,i_2\,i_3\,i_4}=\sum^2_{j_1\!,k_1}
C_{i_1\,j_1\,i_2\,k_1\,i_3\,i_4}\,d^{\kern 1pt k_1\,j_1}.
\mytag{7.5}
$$
The formula \mythetag{7.5} defines a pseudotensorial object of the type
$(0,4)$ and of the weight $4$. Its components can be calculated explicitly
if needed.\par 
     The square shape $\delta$ in Fig\.~7.1 is associated with a pseudoscalar
object of the wight $8$. It is calculated by means of the following
formula:
$$
\hskip -2em
\delta\,\,\,=\!\!\!\sum^2\Sb k_1\!,k_2,k_3,k_4\\j_1\!,j_2,j_3,j_4\endSb
\!\!\!\!\!D_{j_4\,k_1\,j_1\,k_2\,j_2\,k_3\,j_3\,k_4}\,d^{\kern 1pt k_1\,j_1}
\,d^{\kern 1pt k_2,j_2}\,d^{\kern 1pt k_3\,j_3}\,d^{\kern 1pt k_4\,j_4}. 
\mytag{7.6}
$$
It turns our that $\delta$ from \mythetag{7.6} is expressed through $\beta$
from \mythetag{6.6}:
$$
\delta=\frac{\beta^{\kern 1pt 2}}{2}.
$$\par
     Pentagonal shapes cannot be associated with the parameter $I_0$. Therefore
we proceed to hexagonal ones. The shape $E$ in Fig\.~7.1 is presented by the 
formula
$$
\hskip -2em
\gathered
E_{i_1\,i_2\,i_3\,i_4\,i_5\,i_6\,i_7\,i_8\,i_9\,i_{10}\,i_{11}\,i_{12}}
\,\,\,\,\,=\!\!\!\!\!\!\sum^2\Sb k_1,k_2,k_3,k_4,k_5,k_6\\ 
j_1,j_2,j_3,j_4,j_5,j_6\endSb\!\!\!\!
A_{i_1\,i_2\,j_6\,k_1}\,d^{\kern 1pt k_1\,j_1}\,\times\\
\vspace{1ex}
\times\,A_{i_3\,i_4\,j_1\,k_2}\,d^{\kern 1pt k_2,j_2}
A_{i_5\,i_6\,j_2\,k_3}\,d^{\kern 1pt k_3\,j_3}
\,A_{i_7\,i_8\,j_3\,k_4}\,d^{\kern 1pt k4\,j4}\,\times\\
\vspace{1ex}
\times\,A_{i_9\,i_{10}\,j_4\,k_5}\,d^{\kern 1pt k_5,j_5}
A_{i_{11}\,i_{12}\,j_5\,k_6}\,d^{\kern 1pt k_6\,j_6}.
\endgathered
\mytag{7.7}
$$
It defines a pseudotensorial object of the type $(0,12)$ and the
weight $6$. Using \mythetag{7.7}, one can calculate the pseudoscalar 
$\varepsilon_0$ associated with the double hexagon: 
$$
\hskip -2em
\gathered
\varepsilon_0\,\,\,\,\,\,\,=\!\!\!\!\!\!\!\sum^2\Sb k_1\!,k_2,k_3,k_4,k_5,k_6\\
j_1\!,j_2,j_3,j_4,j_5,j_6\endSb
\!\!\!\!\!E_{j_6\,k_1\,j_1\,k_2\,j_2\,k_3\,j_3\,k_4,\,j_4\,k_5\,j_5\,k_6}
\,d^{\kern 1pt k_1\,j_1}\,\times\\
\text{\vphantom{a}}\kern 80pt \times\,d^{\kern 1pt k_2,j_2}\,d^{\kern 1pt k_3\,j_3}
\,d^{\kern 1pt k_4\,j_4}\,d^{\kern 1pt k_5\,j_5}\,d^{\kern 1pt k_6\,j_6}. 
\endgathered
\qquad
\mytag{7.8}
$$
Actually, using \mythetag{7.7} is rather time consuming. For this reason, instead of \mythetag{7.8}, we use another formula for calculating $\varepsilon_0$:
$$
\hskip -2em
\gathered
\varepsilon_0\,\,\,\,\,\,\,=\!\!\!\!\!\!\!\sum^2\Sb k_1\!,k_2,k_3,k_4,k_5,k_6
\\j_1\!,j_2,j_3,j_4,j_5,j_6\endSb
\!\!\!\!\!
B_{j_5\,j_6\,k_1\,k_2}\,B_{j_1\,j_2\,k_3\,k4}\,B_{j_3,j_4,k_5,k_6}
\,d^{\kern 1pt k_1\,j_1}\,\times\\
\text{\vphantom{a}}\kern 80pt \times\,d^{\kern 1pt k_2,j_2}\,d^{\kern 1pt k_3\,j_3}
\,d^{\kern 1pt k_4\,j_4}\,d^{\kern 1pt k_5\,j_5}\,d^{\kern 1pt k_6\,j_6}. 
\endgathered
\qquad
\mytag{7.9}
$$
The quantity \mythetag{7.9} is a pseudoscalar object of the weight $12$. Here is an
explicit expression for this pseudoscalar object:
$$
\gathered
\varepsilon_0=
128\,(A_{1112})^{\kern 1pt 3}\,(A_{1222})^{\kern 1pt 3}
-30\,(A_{1111})^{\kern 1pt 2}\,(A_{2222})^{\kern 1pt 2}
\,(A_{1122})^{\kern 1pt 2}\,-\\
-\,30\,A_{1111}\,(A_{1122})^{\kern 1pt 4}\,A_{2222}
-24\,A_{1111}\,(A_{1122})^{\kern 1pt 3}
\,(A_{1222})^{\kern 1pt 2}\,-\\
-\,12\,(A_{1112})^{\kern 1pt 4}\,(A_{2222})^{\kern 1pt 2}
+96\,A_{1111}\,A_{1112}\,A_{1222}\,(A_{1122})^{\kern 1pt 2}
\,A_{2222}\,+\\
+\,48\,(A_{1112})^{\kern 1pt 3}\,A_{1222}\,A_{1122}\,A_{2222}
+24\,(A_{1111})^{\kern 1pt 2}\,A_{1112}
\,A_{1222}\,(A_{2222})^{\kern 1pt 2}\,+\\
+\,48\,A_{1111}\,A_{1112}\,(A_{1222})^{\kern 1pt 3}\,A_{1122}
+24\,A_{1111}\,A_{2222}^{\kern 1pt 2}\,A_{1122}
\,(A_{1112})^{\kern 1pt 2}\,+\\
+24\,A_{1111}^{\kern 1pt 2}\,A_{2222}\,A_{1122}\,(A_{1222})^{\kern 1pt 2}
-336\,(A_{1122})^{\kern 1pt 2}\,(A_{1112})^{\kern 1pt 2}
\,(A_{1222})^{\kern 1pt 2}\,-\\
-24\,(A_{1122})^{\kern 1pt 3}\,(A_{1112})^{\kern 1pt 2}\,A_{2222}
-120\,A_{1111}\,(A_{1112})^{\kern 1pt 2}\,(A_{1222})^{\kern 1pt 2}
\,A_{2222}\,-\\
-12\,(A_{1111})^{\kern 1pt 2}\,(A_{1222})^{\kern 1pt 4}
-2\,(A_{1111})^{\kern 1pt 3}\,(A_{2222})^{\kern 1pt 3}
-66\,(A_{1122})^{\kern 1pt 6}\,+\\
+264\,A_{1112}\,A_{1222}\,(A_{1122})^{\kern 1pt 4}.
\endgathered\quad
\mytag{7.10}
$$\par
    Formulas similar to \mythetag{7.9} are available for the quantities
from $\varepsilon_1$ through $\varepsilon_{10}$:
$$
\allowdisplaybreaks
\gather
\hskip -2em
\gathered
\varepsilon_1\,\,\,\,\,\,\,=\!\!\!\!\!\!\!\sum^2\Sb k_1\!,k_2,k_3,k_4,k_5,k_6\\
j_1\!,j_2,j_3,j_4,j_5,j_6\endSb
\!\!\!\!\!
C_{j_1\,j_6\,j_2\,j_3\,j_4\,j_5}\,C_{k_1\,k_2\,k_3\,k_4\,k_5\,k_6}
\,d^{\kern 1pt k_1\,j_1}\,\times\\
\text{\vphantom{a}}\kern 80pt \times\,d^{\kern 1pt k_2,j_2}\,d^{\kern 1pt k_3\,j_3}
\,d^{\kern 1pt k_4\,j_4}\,d^{\kern 1pt k_5\,j_5}\,d^{\kern 1pt k_6\,j_6}, 
\endgathered
\qquad
\mytag{7.11}\\
\hskip -2em
\gathered
\varepsilon_2\,\,\,\,\,\,\,=\!\!\!\!\!\!\!\sum^2\Sb k_1\!,k_2,k_3,k_4,k_5,k_6\\
j_1\!,j_2,j_3,j_4,j_5,j_6\endSb
\!\!\!\!\!
\hat C_{j_1\,j_6\,j_3\,j_2\,j_4\,j_5}\,C_{k_1\,k_2\,k_3\,k_4\,k_5\,k_6}
\,d^{\kern 1pt k_1\,j_1}\,\times\\
\text{\vphantom{a}}\kern 80pt \times\,d^{\kern 1pt k_2,j_2}\,d^{\kern 1pt k_3\,j_3}
\,d^{\kern 1pt k_4\,j_4}\,d^{\kern 1pt k_5\,j_5}\,d^{\kern 1pt k_6\,j_6}, 
\endgathered
\qquad
\mytag{7.12}\\
\hskip -2em
\gathered
\varepsilon_3\,\,\,\,\,\,\,=\!\!\!\!\!\!\!\sum^2\Sb k_1\!,k_2,k_3,k_4,k_5,k_6\\
j_1\!,j_2,j_3,j_4,j_5,j_6\endSb
\!\!\!\!\!
\hat C_{j_1\,j_2\,j_3\,j_4\,j_5\,k_1}\,\hat C_{k_6\,k_5\,k_4\,j_6\,k_3\,k_2}
\,d^{\kern 1pt k_1\,j_1}\,\times\\
\text{\vphantom{a}}\kern 80pt \times\,d^{\kern 1pt k_2,j_2}\,d^{\kern 1pt k_3\,j_3}
\,d^{\kern 1pt k_4\,j_4}\,d^{\kern 1pt k_5\,j_5}\,d^{\kern 1pt k_6\,j_6}, 
\endgathered
\qquad
\mytag{7.13}\\
\hskip -2em
\gathered
\varepsilon_4\,\,\,\,\,\,\,=\!\!\!\!\!\!\!\sum^2\Sb k_1\!,k_2,k_3,k_4,k_5,k_6\\
j_1\!,j_2,j_3,j_4,j_5,j_6\endSb
\!\!\!\!\!
\hat C_{j_1,j_2,j_3,j_4,j_5,k_1}\,\hat C_{k_4,k_6,k_3,j_6,k_2,k_5}
\,d^{\kern 1pt k_1\,j_1}\,\times\\
\text{\vphantom{a}}\kern 80pt \times\,d^{\kern 1pt k_2,j_2}\,d^{\kern 1pt k_3\,j_3}
\,d^{\kern 1pt k_4\,j_4}\,d^{\kern 1pt k_5\,j_5}\,d^{\kern 1pt k_6\,j_6}, 
\endgathered
\qquad
\mytag{7.14}\\
\hskip -2em
\gathered
\varepsilon_5\,\,\,\,\,\,\,=\!\!\!\!\!\!\!\sum^2\Sb k_1\!,k_2,k_3,k_4,k_5,k_6\\
j_1\!,j_2,j_3,j_4,j_5,j_6\endSb
\!\!\!\!\!
\hat C_{j_1\,j_3\,j_5\,j_6\,j_4\,k_1}\,\hat C_{k_3\,k_4\,j_2\,k_5\,k_6\,k_2}
\,d^{\kern 1pt k_1\,j_1}\,\times\\
\text{\vphantom{a}}\kern 80pt \times\,d^{\kern 1pt k_2,j_2}\,d^{\kern 1pt k_3\,j_3}
\,d^{\kern 1pt k_4\,j_4}\,d^{\kern 1pt k_5\,j_5}\,d^{\kern 1pt k_6\,j_6}, 
\endgathered
\qquad
\mytag{7.15}\\
\hskip -2em
\gathered
\varepsilon_6\,\,\,\,\,\,\,=\!\!\!\!\!\!\!\sum^2\Sb k_1\!,k_2,k_3,k_4,k_5,k_6\\
j_1\!,j_2,j_3,j_4,j_5,j_6\endSb
\!\!\!\!\!
\hat C_{j_3,j_4,j_1,j_5,j_6,k_1}\,\hat C_{k_3,k_5,j_2,k_4,k_6,k_2}
\,d^{\kern 1pt k_1\,j_1}\,\times\\
\text{\vphantom{a}}\kern 80pt \times\,d^{\kern 1pt k_2,j_2}\,d^{\kern 1pt k_3\,j_3}
\,d^{\kern 1pt k_4\,j_4}\,d^{\kern 1pt k_5\,j_5}\,d^{\kern 1pt k_6\,j_6}, 
\endgathered
\qquad
\mytag{7.16}\\
\hskip -2em
\gathered
\varepsilon_7\,\,\,\,\,\,\,=\!\!\!\!\!\!\!\sum^2\Sb k_1\!,k_2,k_3,k_4,k_5,k_6\\
j_1\!,j_2,j_3,j_4,j_5,j_6\endSb
\!\!\!\!\!
\hat C_{j_1\,j_6\,j_3\,j_4\,j_5\,k_1}\,\hat C_{j_2\,k_6\,k_4\,k_3\,k_5\,k_2}
\,d^{\kern 1pt k_1\,j_1}\,\times\\
\text{\vphantom{a}}\kern 80pt \times\,d^{\kern 1pt k_2,j_2}\,d^{\kern 1pt k_3\,j_3}
\,d^{\kern 1pt k_4\,j_4}\,d^{\kern 1pt k_5\,j_5}\,d^{\kern 1pt k_6\,j_6}, 
\endgathered
\qquad
\mytag{7.17}\\
\hskip -2em
\gathered
\varepsilon_8\,\,\,\,\,\,\,=\!\!\!\!\!\!\!\sum^2\Sb k_1\!,k_2,k_3,k_4,k_5,k_6\\
j_1\!,j_2,j_3,j_4,j_5,j_6\endSb
\!\!\!\!\!
\hat C_{j_1\,j_3\,j_4\,j_5\,j_6\,k_1}\,\hat C_{j_2\,k_3\,k_4\,k_5\,k_6\,k_2}
\,d^{\kern 1pt k_1\,j_1}\,\times\\
\text{\vphantom{a}}\kern 80pt \times\,d^{\kern 1pt k_2,j_2}\,d^{\kern 1pt k_3\,j_3}
\,d^{\kern 1pt k_4\,j_4}\,d^{\kern 1pt k_5\,j_5}\,d^{\kern 1pt k_6\,j_6}, 
\endgathered
\qquad
\mytag{7.18}\\
\hskip -2em
\gathered
\varepsilon_9\,\,\,\,\,\,\,=\!\!\!\!\!\!\!\sum^2\Sb k_1\!,k_2,k_3,k_4,k_5,k_6\\
j_1\!,j_2,j_3,j_4,j_5,j_6\endSb
\!\!\!\!\!
\hat C_{j_1\,j_3\,j_4\,j_5\,j_6\,k_1}\,\hat C_{k_3\,k_4\,j_2\,k_5\,k_6\,k_2}
\,d^{\kern 1pt k_1\,j_1}\,\times\\
\text{\vphantom{a}}\kern 80pt \times\,d^{\kern 1pt k_2,j_2}\,d^{\kern 1pt k_3\,j_3}
\,d^{\kern 1pt k_4\,j_4}\,d^{\kern 1pt k_5\,j_5}\,d^{\kern 1pt k_6\,j_6}, 
\endgathered
\qquad
\mytag{7.19}\\
\hskip -2em
\gathered
\varepsilon_{10}\,\,\,\,\,\,\,=\!\!\!\!\!\!\!\sum^2\Sb k_1\!,k_2,k_3,k_4,k_5,k_6\\
j_1\!,j_2,j_3,j_4,j_5,j_6\endSb
\!\!\!\!\!
\hat C_{j_1\,j_2\,j_3\,j_4\,j_5\,j_6}\,\hat C_{k_1\,k_6\,k_3\,k_2\,k_4\,k_5}
\,d^{\kern 1pt k_1\,j_1}\,\times\\
\text{\vphantom{a}}\kern 80pt \times\,d^{\kern 1pt k_2,j_2}\,d^{\kern 1pt k_3\,j_3}
\,d^{\kern 1pt k_4\,j_4}\,d^{\kern 1pt k_5\,j_5}\,d^{\kern 1pt k_6\,j_6}. 
\endgathered
\qquad
\mytag{7.20}
\endgather
$$
It turns out that $\varepsilon_2$ in \mythetag{7.12} vanishes, i\.\,e\,
we have the equality 
$$
\varepsilon_2=0.
$$ 
Using \mythetag{7.11}, \mythetag{7.13}, \mythetag{7.14}, \mythetag{7.15},
\mythetag{7.16},\mythetag{7.17}, \mythetag{7.18}, \mythetag{7.19},
\mythetag{7.20}, one can derive explicit formulas for $\varepsilon_1$
and for the quantities from $\varepsilon_3$ through $\varepsilon_{10}$.
For instance, in the case of the quantity $\varepsilon_5$ we have
$$
\gathered
\varepsilon_5=
12\,A_{1111}\,(A_{1122})^4\,A_{2222}
-6\,(A_{1111})^{\kern 1pt 2}\,(A_{2222})^{\kern 1pt 2}
\,(A_{1122})^{\kern 1pt 2}\,-\\
-\,12\,A_{1111}\,(A_{1122})^{\kern 1pt 3}\,(A_{1222})^{\kern 1pt 2}
-6\,(A_{1112})^{\kern 1pt 4}\,(A_{2222})^{\kern 1pt 2}\,-\\
-\,24\,A_{1111}\,A_{1112}\,A_{1222}\,(A_{1122})^{\kern 1pt 2}\,A_{2222}
+24\,(A_{1112})^{\kern 1pt 3}\,A_{1222}\,A_{1122}\,A_{2222}\,+\\
+\,24\,A_{1111}\,A_{1112}\,(A_{1222})^{\kern 1pt 3}\,A_{1122}
+12\,A_{1111}\,(A_{2222})^{\kern 1pt 2}\,A_{1122}\,A_{1112})^{\kern 1pt 2}\,+\\
+\,12\,A_{1111})^{\kern 1pt 2}\,A_{2222}\,A_{1122}\,(A_{1222})^{\kern 1pt 2}
-24\,(A_{1122})^{\kern 1pt 2}\,(A_{1112})^{\kern 1pt 2}
\,(A_{1222})^{\kern 1pt 2}\,-\\
-\,12\,(A_{1122})^{\kern 1pt 3}\,(A_{1112})^{\kern 1pt 2}\,A_{2222}
-12\,A_{1111}\,(A_{1112})^{\kern 1pt 2}\,(A_{1222})^{\kern 1pt 2}
\,A_{2222}\,-\\
-\,6\,(A_{1111})^{\kern 1pt 2}\,(A_{1222})^{\kern 1pt 4}-6\,(A_{1122})^{\kern 1pt 6}
+24\,A_{1112}\,A_{1222}\,(A_{1122})^{\kern 1pt 4}.
\endgathered
\quad
\mytag{7.21}
$$
Fortunately there is no need to calculate the rest of the quantities
$\varepsilon_1,\,\ldots,\,\varepsilon_{10}$. They are expressed as 
linear combinations of \mythetag{7.10} and \mythetag{7.21}: 
$$
\xalignat 2
&\hskip -2em
\varepsilon_1=\varepsilon_5,
&&\varepsilon_3=2\,\varepsilon_5,\\
&\hskip -2em
\varepsilon_4=-\varepsilon_5,
&&\varepsilon_6=-\varepsilon_0+3\,\varepsilon_5,\\
\vspace{-1.5ex}
\mytag{7.22}\\
\vspace{-1.5ex}
&\hskip -2em
\varepsilon_7=\varepsilon_5,
&&\varepsilon_8=-2\,\varepsilon_5,\\
&\hskip -2em
\varepsilon_9=2\,\varepsilon_5,
&&\varepsilon_{10}=-\varepsilon_0+4\,\varepsilon_5.
\endxalignat
$$
There is a formula similar to \mythetag{7.22} for the parameter $I_0$:
$$
\hskip -2em
I_0=-\frac{1}{2}\,\varepsilon_0+\frac{11}{2}\,\varepsilon_5. 
\mytag{7.23}
$$
It turns out that $\varepsilon_0$ and $\varepsilon_5$ are expressed
through $\beta^{\kern 1pt 3}$ and $\gamma^{\kern 1pt 2}$:
$$
\xalignat 2
&\hskip -2em
\varepsilon_0=\frac{1}{4}\,\beta^{\kern 1pt 3}
-\frac{1}{3}\,\gamma^{\kern 1pt 2},
&&\varepsilon_5=-\frac{1}{6}\,\gamma^{\kern 1pt 2}.
\mytag{7.24}
\endxalignat
$$
Substituting \mythetag{7.24} into \mythetag{7.23}, we derive 
$$
\hskip -2em
I_0=-\frac{1}{8}\,\beta^{\kern 1pt 3}
-\frac{3}{4}\,\gamma^{\kern 1pt 2},
\mytag{7.25}
$$
where $\beta$ and $\gamma$ are given by \mythetag{6.5} and 
\mythetag{7.4}. The equality \mythetag{7.25} shows that $I_0$ is
a pseudoscalar object of the weight $12$ unlike $I_1$, $I_2$, 
$I_3$, $I_4$, $I_5$, $I_6$, $I_7$, $I_8$, which are just certain 
components of pseudotensorial objects. The equality \mythetag{7.25}
is not surprising since $\beta$ and $\gamma$ are two basic invariants
of a quartic form (see \mycite{9}).\par
\head
8. Conclusions. 
\endhead
     The main result of the present paper consists in revealing the 
tensorial nature of the parameters $I_1$, $I_2$, $I_3$, $I_4$, $I_5$, 
$I_6$, $I_7$, $I_8$ in Dickson's and Rees's positivity test for 
quartic forms. This result is expressed by the formulas \mythetag{6.2},
\mythetag{6.4}, \mythetag{6.9}, and \mythetag{6.10}.
It can be further generalized to trivariate quartic forms and to quartic forms
with greater number of variables. \pagebreak As for the tensorial nature of 
the parameter $I_0$ expressed through the formulas \mythetag{7.25}, 
\mythetag{6.5} and \mythetag{7.4}, it is a classical result known probably 
since Issai Schur, F.~Franklin, J.~J.~Sylvester, and David Hilbert. 

\enddocument
\end